\documentclass[11pt,twoside]{article}
\usepackage{latexsym, amssymb,amsthm,amsmath,epsfig,slashed,bbding,color}
\usepackage[all]{xy}
\usepackage[retainorgcmds]{IEEEtrantools}
\usepackage[pdftex]{hyperref}
\usepackage[margin=10pt,font=small,labelfont=bf]{caption}

\definecolor{darkred}{rgb}{0.5,0,0}
\definecolor{darkgreen}{rgb}{0,0.5,0}
\definecolor{darkblue}{rgb}{0,0,0.5}
\definecolor{darkorange}{rgb}{0.3,0.6,0.2}
\definecolor{darkyellow}{rgb}{0.75,0.75,0.2}

\theoremstyle{plain}
\newtheorem{theorem}{Theorem}[section]
\newtheorem{proposition}[theorem]{Proposition}

\theoremstyle{definition}
\newtheorem{definition}[theorem]{Definition}

\newtheorem{lemma}[theorem]{Lemma}
\newtheorem{corollary}[theorem]{Corollary}

\newtheorem{examplex}[theorem]{Example}
\newenvironment{example}
  {\pushQED{\qed}\examplex}
  {\popQED\endexamplex}

\newcounter{construction}

\newenvironment{renumerate}%
{%
\begin{enumerate}}%
{\end{enumerate}%
}%

\newenvironment{remark}%
{\vskip6pt%
\noindent%
{\it Remark.}}%
{\vskip6pt}

\newenvironment{remarks}%
{\vskip6pt%
\noindent%
{\it Remarks}. %
\begin{renumerate}}%
{\end{renumerate}\vskip6pt}

\newenvironment{buildingblock}%
{\vskip6pt%
\noindent%
{\bf Building block:}}%
{\vskip6pt}

\newenvironment{ingredients}%
{\vskip6pt%
\noindent%
{\bf Ingredients:}}%
{\vskip6pt}

\newenvironment{thesurgery}%
{\vskip6pt%
\noindent%
{\bf The surgery:}}%
{\vskip6pt}

\newenvironment{introtheorem}[1]%
{\vskip6pt%
\noindent%
{\bf Theorem  {#1}.}\it}%
{\vskip6pt}

\makeatletter
\def\Ddots{\mathinner{\mkern1mu\raise\p@
\vbox{\kern7\p@\hbox{.}}\mkern2mu
\raise4\p@\hbox{.}\mkern2mu\raise7\p@\hbox{.}\mkern1mu}}
\makeatother

\newcommand{\R}{\text{${\mathbb R}$}}
\newcommand{\C}{\text{$\mathbb C$}}

\newcommand{\Z}{\text{$\mathbb Z$}}

\newcommand{\Q}{\text{$\mathbb Q$}}

\renewcommand{\frak}[1]{\text{$\mathfrak{#1}$}}
\renewcommand{\tilde}{\widetilde}

\newcommand{\N}{\text{$\mathcal{N}$}}

\newcommand{\U}{\text{$\mathcal{U}$}}

\newcommand{\e}{\text{$\varepsilon$}}

\renewcommand{\bar}{\overline}

\newcommand{\gf}{\text{$\varphi$}}

\newcommand{\Id}{\mathrm{Id}}

\newcommand{\del}{\text{$\partial$}}

\newcommand{\tensor}{\otimes}

\newcommand{\mc}[1]{\text{$\mathcal{#1}$}}

\newcommand{\into}{\longrightarrow}
\newcommand{\noqed}{\let\qed\relax}

\newcommand{\IP}[1]{\langle #1 \rangle}

\newcommand{\nhood}{neighbourhood}
\newcommand{\nhoods}{neighbourhoods}

\newcommand{\wrt}{with respect to}
\renewcommand{\iff}{if and only if}

\date{} \usepackage{color} \definecolor{tocolor}{rgb}{.1,.1,.5}
\definecolor{urlcolor}{rgb}{.2,.2,.6}
\definecolor{linkcolor}{rgb}{.1,.1,.6}
\definecolor{citecolor}{rgb}{.6,.2,.1}
\hypersetup{backref=true, colorlinks=true, urlcolor=urlcolor,
  linkcolor=linkcolor, citecolor=citecolor}

\begingroup
\hypersetup{linkcolor=blue}
\endgroup

\numberwithin{equation}{section}

\begin{document}

\title{Examples and counter-examples of  log-symplectic manifolds}
\author{
Gil R. Cavalcanti\thanks{{\tt gil.cavalcanti@gmail.com}} \\[6pt]
Department of Mathematics\\
Utrecht University\\
}
\maketitle

\abstract{We study topological properties of log-symplectic structures and produce examples of compact manifolds with such structures. Notably we show that several symplectic manifolds do not admit {\it bona fide} log-symplectic structures and several {\it bona fide} log-symplectic manifolds do not admit symplectic structures, for example $\#m \C P^2 \# n\bar{\C P^2}$ has {\it bona fide} log-symplectic structures \iff\ $m,n>0$ while they only have symplectic structures for $m=1$. We introduce surgeries that produce log-symplectic manifolds out of symplectic manifolds  and show that any compact oriented log-symplectic four-manifold can be transformed into a collection of symplectic manifolds by reversing these surgeries. Finally we show that if a compact manifold admits an achiral Lefschetz fibration with homologicaly essential fibers, then the manifold admits a log-symplectic structure. Then, using  results of Etnyre and Fuller \cite{MR2219214}, we conclude that if $M$ is  a compact, simply connected 4-manifold then $M \# (S^2 \times S^2)$ and $M \# \C P^2 \# \bar{\C P}^2$ have log-symplectic structures.}


\tableofcontents

\section{Introduction}

A log-symplectic structure on a manifold $M^{2n}$ is a Poisson structure $\pi \in \frak{X}^2(M)$ for which $\pi^n$ has only nondegenerate zeros. This condition is weaker than asking that $M$ is outright symplectic (in which case $\pi^n$ would not vanish) and yet it is only a little less so, since it still requires that $\pi$ is generically symplectic and that  its failure to be so everywhere is as well behaved as one could ask. If we want to rule out log-symplectic structures which are in fact symplectic, we refer to them as {\it bona fide} or nonsymplectic log-symplectic structures.

These structures have been classified on surfaces by  Radko  \cite{MR1959058}  and already in dimension two there is a marked contrast with the symplectic case, namely, every surface (orientable or not) has a log-symplectic structure.  Recently log-symplectic structures received renewed attention: Guillemin, Miranda and Pires  \cite{Guilleminetal2012} proved a local form for the Poisson structure in a \nhood\ of the zeros of $\pi^n$ and Gualtieri and Li  \cite{Gualtieri-Li-2012} managed to give a clear geometrical description of symplectic groupoids integrating log-symplectic structures.

Despite these recent advances in the theory, the area still lacks examples and even topological obstructions to the existence of these structures are unknown. So, given a manifold, the question ``does it have a log-symplectic structure?" is a little hard to answer. 

We tackle these shortcomings in this paper. Indeed, Marcut and Osorno-Torres's paper \cite{marcut-osorno,osorno-thesis} and the present one are the first to provide topological obstructions to the existence of log-symplectic structures. While Marcut and Osorno-Torres prove that a log-symplectic manifold whose singular locus has a compact component must have a cohomology class $a \in H^2(M)$ such that $a^{n-1}\neq 0$, we prove a different property which is more contrastive with symplectic geometry:

\begin{introtheorem}{\ref{theo:simple invariant}}
If a compact oriented manifold $M^{2n}$, with $n >1$, admits a {\it bona fide} log-symplectic structure then there are classes  $a,b\in H^2(M;\R)$ such that $a^{n-1}b \neq0$ and  $b^2=0$.
\end{introtheorem}
Different from Marcut and Osorno-Torres's  topological constraint, the existence of the class $b$  is not necessarily shared by symplectic manifolds and, in effect, shows that there are several symplectic manifolds for which the only log-symplectic structures are outright symplectic while other manifolds do not admit log-symplectic structures at all.

We then move on to produce examples of manifolds admitting such structures. The first aproach consists simply of deforming a symplectic structure into a {\it bone fide} log-symplectic. We show:
\begin{introtheorem}{\ref{theo:birth}}
Let $(M^{2n},\omega)$ be a symplectic manifold and $k>0$ be an integer. If $M$ has a compact symplectic submanifold $F^{2n-2}\subset M$ with trivial normal bundle, then $M$ has a log-symplectic structure  for which the zero locus of  $\pi^n$ has $k$ components all diffeomorphic to $F\times S^1$.  
\end{introtheorem}

Using symplectic blow-up we can then construct log-symplectic structures on $\#m\C P^2 \# n\bar{\C P^2}$ for $m,n>0$. Therefore coupling the two theorems we have a complete classification of which manifolds in the family $\#m\C P^2 \# n\bar{\C P^2}$ for $m,n\geq0$  admit log-symplectic structures (see Figure \ref{fig:table}).

\begin{figure}[h!!]
\begin{center}
\begin{tabular}{|c|c|c|}
\hline\hline
$\#m \C P^2 \# n \bar{\C P^2}$ &symplectic& {\it bona fide} log-symplectic\\
\hline\hline
$m>1$, $n>0$& {\color{darkred} \XSolidBrush} & {\color{darkgreen} \checkmark}\\
\hline
$m>1$, $n=0$& {\color{darkred} \XSolidBrush}& {\color{darkred} \XSolidBrush}\\
\hline
$m=1$, $n>0$&{\color{darkgreen} \checkmark} &  {\color{darkgreen} \checkmark}\\
\hline
$m=1$, $n=0$& {\color{darkgreen} \checkmark}&  {\color{darkred} \XSolidBrush}\\
\hline
$m=0$, $n>0$& {\color{darkred} \XSolidBrush} & {\color{darkred} \XSolidBrush}\\
\hline\hline
\end{tabular}
\end{center}
\caption{Table showing the values of $m$ and $n$ for which $m \C P^2 \# n \bar{\C P^2}$ has symplectic or {\it bona fide} log-symplectic structures. In the symplectic case, we require that the orientation determined by the symplectic structure agrees with the orientation of the manifold. In the log-symplectic case, since these structures do not induce a preferred orientation on the manifold, we simply assert the existence of the structure in the underlying unoriented differentiable manifold.}\label{fig:table}
\end{figure}

Further we introduce two surgeries which produce log-symplectic manifolds out of log-symplectic manifolds and which increase the number of components of the singular locus of the Poisson structure, hence even if the starting  manifolds are symplectic, the resulting manifolds will only be log-symplectic.

Following the lines of Gompf's theorem relating symplectic structures to Lefschetz fibrations, we prove an analogue result for log-symplectic manifolds:
\begin{introtheorem}{\ref{theo:achiral}}
Let $M^4$ and $\Sigma^2$ be compact connected manifolds and  $p:M \into \Sigma$ be an achiral Lefschetz fibration with generic fiber $F$. If $F$ is orientable and $[F]\neq 0 \in H_2(M;\R)$, then $M$ has a log-symplectic structure whose singular locus has one component and for which the fibers are symplectic submanifolds of the symplectic leaves of the Poisson structure.
\end{introtheorem}

Using results of Etnyre and Fuller on such fibrations \cite{MR2219214} we obtain a general existence result

\begin{introtheorem}{\ref{theo:corollary from Etnyre--Fuller}}
Let $M$ be a simply connected compact four-manifold. Then both $M \# (S^2 \times S^2)$ and $M\#\C P^2 \# \bar{\C P^2}$ admit bona fide log-symplectic structures. 
\end{introtheorem}

We finish showing that in four dimensions any compact orientable log-symplectic manifold is obtained out of a symplectic manifold using our surgeries. Expressed another way:

\begin{introtheorem}{\ref{theo:4d}}
Let $(M^4,\pi)$ be a compact, orientable, log-symplectic manifold with singular locus $Z$. Then each unoriented component of $M\backslash Z$ is symplectomorphic to an open subset of  a compact symplectic manifold
\end{introtheorem}

While our original motivation to study log-symplectic structures lies in the realm of Poisson geometry,   one might naturally bundle them together with folded symplectic structures: other structures defined as degenerate symplectic structures whose degeneracy locus is defined by a transversality condition. In particular it is natural to compare these structures and therefore put our results in context, specially because there are a few results similar in content to the ones we obtain here.

The relevant result regarding  {\it topological obstructions} of folded symplectic structures  was proved by Cannas da Silva in  \cite{MR2670164}: every compact oriented 4-manifold admits a folded symplectic structure, hence, differently from log-symplectic structures, in 4-dimensions there are no topological obstructions to the existence of folded symplectic structures. The result on {\it symplectization} of folded symplectic structures was proved in \cite{MR1748286} by Cannas da Silva {\it et al}. Differently from log-symplectic structures (c.f. Theorem \ref{theo:4d}),  not all folded symplectic structures can be   ``unfolded" as a condition on the one dimensional foliation of the folding must be imposed.  Finally, a {\it relation between folded symplectic structures and achiral Lefschetz fibrations} was obtained by Baykur in \cite{MR2253445}: achiral Lefschetz fibrations with homologically nontrivial fibers admit folded symplectic structures compatible with the fibration. This result is analogous in statement and proof to our Theorem \ref{theo:achiral}, as both proofs are based on Gompf's original result for symplectic manifolds.  Notice that our result is stronger than Baykur's because, in any dimension, if a manifold admits a log-symplectic structure, it also admits a folded symplectic structure whose folding is the singular locus of the log-symplectic structure \cite{Guilleminetal2015}.

This paper is organized as follows.  Section \ref{sec:basics} reviews the basics of Poisson geometry relevant for our study and Section \ref{sec:localform} reviews Guillemin--Miranda--Pires normal form theorem \cite{Guilleminetal2012}. Section \ref{sec:simpleinvariant} introduces a simple topological invariant that allows us to show that there are many symplectic manifolds which do not admit {\it bona fide} log-symplectic structures. Section \ref{sec:birth} shows that under general assumptions one can deform a symplectic structure into a log-symplectic structure and Section \ref{sec:surgeries} introduces the surgeries and gives the existence result for log-symplectic structures on achiral Lefschetz fibrations. Finally, Section \ref{sec:reversing} shows that in four dimensions the surgeries can be reversed and any compact, orientable, log-symplectic four manifold can be transformed into a symplectic manifold by surgeries.

\noindent
{\bf Acknowledgements}: The author is thankful to Ioan Marcut, Ralph Klaasse, Pedro Frejlich, Marius Crainic for useful conversations. The author is specially thankful to Ioan Marcut for the argument of Theorem \ref{theo:proper} and for explaining the results from \cite{marcut-osorno}.

While carrying out this research, the author was made aware that Frejlich, Martinez-Torres and Miranda were carrying out a project \cite{FMM13} which overlaps with the results in this paper. Notably, they had independently produced our ``Construction \ref{cons:singular fiber}" and our Theorem \ref{theo:4d}.

This research is supported by a VIDI grant from the Dutch Science Foundation (NWO).

\section{Poisson structures}\label{sec:basics}

This section we give a short account of the basic material on Poisson and log-symplectic structures. For more details we refer the reader to \cite{Gualtieri-Li-2012,MR2861781,Guilleminetal2012}.

\subsection{Poisson cohomologies}

A Poisson structure on a manifold $M^m$ is a bivector $\pi \in \frak{X}^2(M) = \Gamma(\wedge^2TM)$ for which
$$[\pi,\pi]=0,$$
where the bracket used is the Schouten--Nijenhuis bracket of multivector fields. Assuming that $M^m$ is even dimensional, say, $m=2n$, a generic bivector (at a point) would give an isomorphism $\pi:T_p^*M \stackrel{\cong}{\into} T_pM$. In this case $\pi^n$ is a non zero element in $\wedge^{2n}T_pM$. If a Poisson bivector $\pi$ is {\it everywhere generic} (i.e., everywhere invertible) then the 2-form $\omega = \pi^{-1}$ is a symplectic structure on $M$. 

\begin{definition}
 The locus where $\pi:T^*M \into TM$ is an isomorphism is the {\it symplectic locus} and  its complement is the  {\it singular locus} of the Poisson structure.
 \end{definition}

Moving to a move general situation which is still modeled on a ``generic" case, one can require that $\pi^n$ only has nondegenerate zeros:

\begin{definition}
 A  {\it log-symplectic structure} on $M^{2n}$ is a Poisson structure $\pi$ for which the zeros of $\pi^n$ are nondegenerate.
\end{definition}

It follows from the definition, using Weinstein's splitting theorem, that one can find coordinates in a \nhood\ of any singular point which render a log-symplectic structure $\pi$ in the following form
$$\pi = x_1 \del_{x_1} \wedge \del_{x_2} + \del_{x_3}\wedge \del_{x_4} + \cdots \del_{x_{2n-1}} \wedge \del_{x_{2n}},$$
and its inverse is given by
$$ \omega = d\log |x_1| \wedge dx_2 + dx_3 \wedge dx_4 + \cdots dx_{2n-1} \wedge dx_{2n}.$$
The fact that the ``symplectic form" $\omega$ acquires a logarithmic singularity along the singular locus of $\pi$ justifies the name of the structure.


Continuing with the general theory, any Poisson manifold comes equipped with two diferential operators which give rise to cohomology theories. The first is the  Poisson differential on multivector fields:
$$d_\pi:\mathfrak{X}^{\bullet}(M)  \into  \mathfrak{X}^{\bullet+1}(M);\qquad d_\pi(\xi) = [\pi,\xi],$$
The Poisson condition and the Jacobi identity for the Schouten--Nijenhuis bracket imply that $d_\pi^2=0$ and its cohomology is known as the {\it Poisson cohomology} of $(M,\pi)$.

The second is the Koszul differential on forms:
$$\delta:\Omega^{\bullet}(M)  \into  \Omega^{\bullet-1}(M);\qquad \delta\rho = \{\pi,d\} \rho,$$
where $\{\pi,d\} = \pi d  - d \pi$ is the graded commutator of operators and $\pi$ acts on forms by interior product. Again the Jacobi identity for the graded commutator and the Poisson condition imply that $\delta^2=0$ and its cohomology is known as the {\it canonical cohomology} of $(M,\pi)$.

These operators are related:
\begin{lemma}
Let $(M,\pi)$ be a Poisson manifold, $\xi \in  \mathfrak{X}^{\bullet}(M)$ and $\rho \in \Omega^{\bullet}(M)$. Then
$$\{\delta,\xi\} \rho = (d_\pi \xi)\cdot \rho,$$
where $\xi$ and $d_\pi \xi$ act on forms by inner product.
\end{lemma}
\begin{proof}
This follows automatically from the description of the Schouten--Nijnehuis bracket as a derived bracket:
$$\{\delta,\xi\} \rho = \{\{\pi,d\},\xi\}\rho  = [\pi,\xi]\cdot \rho = (d_\pi \xi)\rho.$$
\end{proof}

\subsection{The canonical bundle and the modular vector field}

Given a Poisson manifold $(M^m,\pi)$, the determinant bundle $K = \wedge^m T^*M$ is also known as the {\it canonical bundle} of $M$. Given any nonvanishing local section $\rho \in \Gamma(K)$ there is a unique vector field $X$ such that
$$\delta \rho = \iota_X \rho.$$
The vector field  $X$ is called the {\it modular vector field}. Notice that changing the trivialization $\rho$ by a nonvanishing function, say $g$,  changes the modular vector field from $X$ to $X + \pi(d\log|g|) = X + d_\pi \log|g|$. In particular, changing  $\rho$ to $-\rho$ does not change $X$ and a modular vector field is determined by a section of the quotient sheaf $K/\Z_2$. If $M$ is nonorientable, there is no global nonvanishing section of $K$, yet, $K/\Z_2$, the sheaf of densities, always has a nonvanishing section, so one can always find globally defined modular vector fields. 

Notice that any modular vector field $X$ is an infinitesimal symmetry of the Poisson structure, i.e., $[\pi,X] =0$  since for a local section $\rho \in \Gamma(K_\pi)$ we have
$$ 0 = \delta^2 \rho = \delta(X \cdot \rho) = (d_{\pi}X)\cdot \rho - X \cdot(X\cdot \rho) = (d_\pi X)\cdot \rho,$$
which implies $[\pi,X]=0$. An immediate consequence is that the rank of the Poisson structure along the flow of a point is constant. Since a different choice of section of $K_\pi/\Z_2$ changes $X$ to $X + d_\pi log|g|$ we see that the modular vector field gives a well defined degree one Poisson cohomology class. A Poisson structure is {\it unimodular} if this class is trivial, which is therefore equivalent to the existence of a globally defined $\delta$-closed section of $K_\pi/\Z_2$.

\begin{definition}
A {\it representation} of a Poisson structure is a vector bundle $E \into (M,\pi)$ together with a flat Poisson connection $\nabla:\Gamma(E)\into \Gamma(TM \tensor E)$, i.e., for $f\in C^\infty(M)$ and $s \in \Gamma(E)$
$$\nabla (f s) = d_\pi f s + f \nabla s\qquad \mbox{ and }\qquad \nabla^2 =0.$$
\end{definition}

\begin{example}
The canonical bundle of a Poisson manifold is a representation. Indeed, the operator $\delta:K \into \Omega^{m-1}(M) \cong \Gamma(TM\tensor K)$ satisfies the properties required for a connection and $\delta^2 =0$ is the flatness condition. Note that if $M$ is orientable, a Poisson structure is unimodular \iff\ its canonical bundle is the trivial representation.
\end{example}

\subsection{Log-symplectic structures --- basics}

Now we can focus on the objects in which we are interested.  Since in a log-symplectic manifold the singular locus is given by the nondegenerate zero locus of a section of a line bundle, we have that it is a smooth submanifold of codimension one. Further, as the rank of the Poisson structure does not change along each of its symplectic leaves, we see that each connected component of the singular locus is itself a Poisson submanifold of $M$, i.e., a union of symplectic leaves.

The following proposition adds up the basic facts about the singular locus.
\begin{proposition}\label{prop:oriented}
Let $M$ be a  log-symplectic manifold, $Z$ its singular locus, $\N^*_Z$ the conormal bundle of $Z$ and $K_M$ and $K_Z$ the canonical bundles of $M$ and $Z$, respectively. Then
\begin{enumerate}
\item $Z$ is an orientable Poisson submanifold of $M$ with symplectic leaves of codimension 1;
\item  $K_Z$  is the trivial representation and has a distinguished trivialization;
\item  $K_M|_Z$ is a Poisson representation over $Z$.
\item  $\mc{N}^*_Z \cong K_M|_Z$ as vector bundles and hence $\mc{N}^*_Z$  inherits the structure of a Poisson representation 
\end{enumerate}
In particular if $M$ is orientable each component of $Z$ has trivial normal bundle. 
\end{proposition}
\begin{proof}
 We have already argued most of the claim 1). The rest follows from the normal form for a log-symplectic structure. Indeed, if 
 $$ \pi = x_1\del_{x_1} \wedge \del_{x_2} + \del_{x_3}\wedge \del_{x_4} + \cdots + \del_{x_{2n-1}}\wedge \del_{x_{2n}} $$ 
then the induced Poisson structure on the singular locus, $[x_1=0]$, is 
 $$\del_{x_3}\wedge \del_{x_4} + \cdots + \del_{x_{2n-1}}\wedge \del_{x_{2n}}$$
 which has codimension one leaves.

To prove 2), we let $\omega = \pi^{-1}$. Then $\omega$ is a 2-form with a logarithmic singularity along $Z$ and the desired volume form on $Z$ is just the residue of $\omega^2$ over $Z$.  In the coordinates used above we have
$$\mathrm{Re}(\omega^n) = dx_2 \wedge\cdots\wedge d x_{2n},$$
hence the residue of $\omega^n$ over $Z$ is nowhere vanishing and one can readily compute $\delta_Z \mathrm{Res}\omega =0$, in these coordinates. 
 
  Claim 3) follows from the fact that a Poisson representation $(E,\nabla)$ over $M$ induces a representation on a Poisson submanifold $Z$ \iff\ for every local section $\rho \in \Gamma(E)$,  we have $\nabla \rho = X_i \rho_i$, where $\rho_i \in \Gamma(E)$ is a local basis for $E$ and $X_i$ is tangent to $Z$ at all points of $Z$.  In our case, the representation is the canonical bundle, $\rho$ is a local nonvanishing volume form  and $\delta \rho = X \rho$, for $X$ the modular vector field, which is tangent to $Z$ as the rank of the Poisson structure must remain constant along the integral curves of $X$. So claim 3) follows.
 
 As for 4),  since $Z$ is the nondegenerate zero locus of $\pi^n \in \Gamma(\wedge^{2n}TM)$, we have that, over $Z$, $d\pi^n$ gives an isomorphism of vector bundles
$$d\pi^n:\wedge^{2n}TM|_Z \tensor \mc{N}^*_Z\stackrel{\cong}{\into} \R$$
that is $K_M|_Z$ is isomorphic to $\mc{N}^*_Z$.
\end{proof}

\section{Invariants and local forms}\label{sec:localform}

While Proposition \ref{prop:oriented} gives a list of simple invariants associated to a log-symplectic structure in \cite{Guilleminetal2012} Guillemin, Miranda and Pires showed that these are in fact all invariants associated to a \nhood\ of the singular locus. Indeed, the following is a direct consequence of the results in \cite{Guilleminetal2012}:

\begin{theorem}
Let $(M,\pi)$ be a log-symplectic manifold and let $Z$ be a compact connected component of the singular locus. Then a \nhood\ of $Z$ is determined by the Poisson structure induced on $Z$, a distiguished flat section of $K_Z$ and its representation on the conormal bundle of $Z$.
\end{theorem}

Taking the inverse of the Poisson structure, one can translate this information into differential forms (c.f. \cite{Guilleminetal2012}):

\begin{theorem}\label{theo:general local form}
Let $(M,\pi)$ be a log-symplectic manifold, let $Z$ be a connected component of the singular locus and $X$ a modular vector field of $\pi$. Then the pair $(\pi,X)$ determines the following structure on $Z$:
\begin{enumerate}
\item The normal bundle of $Z$ as a vector bundle, i.e., a class $w_1 \in H^1(Z,\Z_2)$.
\item A closed 1-form $\theta \in \Omega^1(Z)$ such that $\theta(X) =-1$.
\item A closed 2-form $\sigma \in \Omega^2(Z)$ such  that $\iota_X \sigma =0$ and
\begin{equation}\label{eq:cosymplectic}
\theta \wedge\sigma^{n-1} \neq 0.
\end{equation}
\end{enumerate}
Changing the modular vector field  by $d_\pi f$ does not change $\theta$ and changes $\sigma$ to $\sigma + df \wedge \theta$.

Further, if $Z$ is compact, any log-symplectic structure inducing the data above on $Z$ is equivalent to a \nhood\ of the zero section of the normal bundle of  $Z$ endowed with the following structure
\begin{equation}\label{eq:general local form}
d\log|x| \wedge \theta + \sigma,
\end{equation}
where  $|\cdot|$ is the distance to the zero section measured \wrt\ a fixed fiberwise linear metric on  $\mc{N}_Z$.
\end{theorem}

Under the conditions of the theorem, the annihilator of the form $\theta$ corresponds to the distribution in $Z$ determined by the Poisson structure and $\sigma$ agrees with the leafwise symplectic form on $Z$.

\begin{definition}
A {\it cosymplectic structure} on a manifold $Z^{2n-1}$  is a pair of closed forms $\theta\in \Omega^1(Z)$ and $\sigma \in \Omega^2(Z)$ satisfying \eqref{eq:cosymplectic}.
\end{definition}

For special types of log-symplectic structure one can rephrase the data {\it 1. --- 3.} above as  a more workable set.
\begin{definition}\label{def:proper}
A connected cosymplectic manifold $(Z,\sigma,\theta)$ is {\it proper} if it is compact and the distribution given by the annihilator of $\theta$ has a compact leaf. A component  $Z$ of the singular locus of a log-symplectic manifold is {\it proper} if the cosymplectic structure induced  on $Z$ is proper. A log-symplectic manifold is {\it proper} if all components of the singular locus are proper. 
\end{definition}
Given a cosymplectic manifold $(Z,\sigma,\theta)$, if we let $X$ be a vector field such that  $\theta(X) = -1$ and $\iota_X \sigma =0$, we have that $\mc{L}_X\theta = 0$ and hence the flow of $X$ preserves the leaves of the distribution determined by $\theta$, hence, if $Z$ is proper with compact (symplectic) leaf  $(F,\sigma) \subset Z$, the flow of $F$ by the vector field $X$ will provide  further leaves of $\pi$. Since $X$ is transverse to $F$ and $Z$ is compact, we see that after finite time, say  $\lambda >0$, the flow  of $X$ brings $F$ back to itself:
$$\gf_{\lambda}:F \into F,$$
Since $\mc{L}_X \sigma =0$, the flow is a symplectomorphism of $F$ and hence $Z$ is a symplectic fiber bundle with fiber  $(F,\sigma)$ over the circle:
$$ Z = \R \times F/\Z,$$
where the quotient is taken \wrt\  the $\Z$-action generated by $(y,p) \mapsto  (y + \lambda,\gf_\lambda(p)).$ Further, the modular vector field is  $-\del_y$ hence $\theta = dy$.

Different choices of nonvanishing sections of $K_\pi/\Z_2$, change the modular vector field over $Z$ by adding Hamitonian vector fields of $F$ so the symplectomorphism $\gf_\lambda$ is only determined up to Hamiltonian symmetries, i.e., the relevant data is only  its class in $\mathrm{Symp}(F)/\mathrm{Ham}(F)$.

Finally, the normal bundle of $Z$ is determined by its first Stiefel--Whitney class $w_1 \in H^1(Z,\Z_2) = H^1(F;\Z_2)^{\gf_\lambda} \times H^1(S^1;\Z_2)$. So, in the proper case, Theorem \ref{theo:general local form} becomes (c.f. \cite{Guilleminetal2012,Gualtieri-Li-2012})

\begin{theorem}\label{theo:proper local form}
Let $Z$ be a proper component of the singular locus of a log-symplectic structure $\pi$ and $F\subset Z$ be a compact symplectic leaf of $\pi$. Then $\pi$ determines the following data:
\begin{enumerate}
\item The normal bundle of $Z$, i.e., a class $w_1 \in H^1(Z,\Z_2) = H^1(F;\Z_2)^{\gf} \times H^1(S^1;\Z_2)$.
\item  The symplectic structure $\sigma$ of $F$;
\item A class $[\gf] \in \mathrm{Symp}(F)/\mathrm{Ham}(F)$; 
\item A period $\lambda >0$;
\end{enumerate}
Further, any two log-symplectic structures inducing the same set of data are equivalent and given a set of data 1. --- 4. there is a proper log-symplectic structure which realises it.
\end{theorem}

Notice that given a nonorientable Poisson manifold, $M$, one can always pass to the oriented double cover $\tilde{M}$ of $M$ which  inherits a Poisson structure from $M$. For the log-symplectic case, this allows us to get a simpler local model for the singular locus as now its \nhood\ depends on one less parameter, since according to Proposition \ref{prop:oriented}   $w_1 =0$ in $\tilde{M}$ .

The following theorem, communicated to the author by Ioan Marcut (see also \cite{osorno-thesis}), uses a Tischler type argument to show one can always deform  a log-symplectic structure into a proper one. In its cosymplectic version,  it had already appeared in \cite{MR2481690}.


\begin{theorem}\label{theo:proper}
If the components of the singular locus of a log-symplectic structure are compact, then the structure can be deformed into a proper one.
\end{theorem}
\begin{proof}
Let $Z$ be a connected component of the singular locus. The proof consist of two steps. Firstly we notice that one can deform the cosymplectic structure $(\theta,\sigma)$ of $Z$ into $(\tilde\theta,\sigma)$ so that the kernel of $\tilde\theta$ gives a fibration structure to $Z$. The second step is to show that this deformation can be realised as a deformation of the log-symplectic structure.

For the first step, let $\tilde{\theta}$ be a closed 1-form representing a class in $H^1(Z,\Q)$  which is close enough to $\theta$ so that we still have
$$\tilde{\theta}\wedge \sigma^{n-1} \neq 0.$$
Since $[\tilde\theta]$ represents a rational class, $[\tilde\theta](H_1(Z;\Z))$ is a lattice $\Lambda$ in $\R$.
Then we define the projection map
$$p:Z \into \R/\Lambda;\qquad p(z) = \int_{z_0}^z \tilde\theta ,$$
where $z_0 \in Z$ is a fixed reference point and the value of the integral modulo $\Lambda$ does not depend of choice of  path  connecting $z_0$ to $z$. By construction  $dp = \tilde\theta$ is nowhere vanishing and hence $p:Z \into S^1$ is a  fibration.

For the second step, according to Theorem \ref{theo:general local form}  there is $\delta>0$ such that the log-symplectic structure in a \nhood\ of $Z$ is equivalent to \eqref{eq:general local form} for $|x|<\delta$. If we let $\psi$ be a smooth function such that 
 $$\psi:[0,1]\into [0,1];\qquad 
\begin{cases}
\psi(x) = 1& \mbox{ if }x < \delta/3\\
\psi(x) = 0& \mbox{ if }x > 2\delta/3\\
\end{cases} $$
then the log-symplectic form
$$d\log|x| \wedge ((1-\psi(|x|)\theta + \psi(|x|)\tilde\theta) +   \sigma,$$
induces the cosymplectic structure $(\tilde\theta,\sigma)$ on $Z$ and agrees with the original log-symplectic structure if $|x|>2\delta/3$ hence can be extended to the rest of $M$ by the original log-symplectic structure.

\end{proof}

\section{A simple topological invariant}\label{sec:simpleinvariant}

One of the simplest and yet restrictive topological properties of  compact symplectic manifolds is the existence of a class $a \in H^2(M)$ whose top power is nonzero. Of course, this does not hold on all log-symplectic manifolds, yet log-symplectic manifolds are just a little shy of satisfying this property as shown by Marcut and Osorno-Torres:

\begin{theorem}[Marcut--Osorno-Torres \cite{marcut-osorno,osorno-thesis}]\label{marcut-Osorno-Torres}
Let $M^{2n}$ be a log-symplectic manifold whose singular locus has a compact component. Then there is a cohomology class $a\in H^2(M;\R)$ such that $a^{n-1} \neq 0$. Further, if $Z \subset M$ is a proper component of the singular locus and has $(F,\sigma)$ as a symplectic leaf, $a$ can be chosen so that  $[a]|_F = [\sigma]$.
\end{theorem}

Here we use a little more of the log-symplectic structure in the orientable case to find another topological property of these manifolds.

\begin{theorem}\label{theo:simple invariant}
If a compact oriented manifold $M^{2n}$, with $n >1$, admits a {\it bona fide} log-symplectic structure then there are classes  $a,b\in H^2(M;\R)$ such that $a^{n-1}b \neq 0$ and  $b^2=0$.
\end{theorem}
\begin{proof}
Assume that $M$ has a log-symplectic structure with singular locus $Z \neq \emptyset$. Then, due to Theorem \ref{theo:proper}, we may assume that the structure is proper, hence $Z$ is a symplectic fibration over the circle with fiber a symplectic manifold $F$. On the one hand, due to Marcut--Osorno-Torres's Theorem there is a globally defined closed 2-form $\tilde\omega \in \Omega^2(M)$ which restricts to the symplectic form on $F$, i.e., the homology class of $F$ pairs nonzero with $a^{n-1}$, so we have that $[F] \neq 0 \in H_{2n-2}(M;\R)$. On the other hand, since, even within $Z$, $F$ appears as a fiber of a fibration, we conclude that the Poincar\'e dual of $F$, $b \in H^2(M;\R)$, must satisfy $b^2 =0$ and, by definition of Poincar\'e dual, 
$$\IP{a^{n-1}b,[M]} = \IP{a^{n-1},F} \neq 0.$$
\end{proof}

A few immediate corollaries:
\begin{corollary}
An orientable, compact, {\it bona fide} log-symplectic manifold $M$ of dimension $2n$ has  $b_{2i}(M) \geq 2$ for $0 < i <n$.
\end{corollary}
\begin{proof}
It follows directly from the relations $a^{n-1} b \neq 0$ and $b^2=0$  that the classes $a^i$ and $a^{i-1}b$ are linearly independent for $0 < i < n$.\end{proof}

\begin{corollary}
For $n>1$, $\C P^n$ has no {\it bona fide} log-symplectic structure and, for $n>2$, the blow-up of $\C P^n$ along a symplectic submanifold of real codimension greater than 4 also does not carry {\it bona fide} log-symplectic structures.
\end{corollary}

\begin{corollary}
A smooth orientable compact four-manifold with definite intersection form  does not admit {\it bona fide} log-symplectic structures. In particular, for $n> 0$, $\#n \C P^2$ and $\#n \bar{\C P^2}$ do not admit {\it bona fide} log-symplectic structures
\end{corollary}
\begin{proof}
Indeed, under the hypothesis of both corollaries there is no element in second cohomology whose square is zero.
\end{proof}

Notice that due to  Taubes result on Seiberg--Witten invariants of symplectic manifolds \cite{MR1306023}, $\#n \C P^2$ does not admit symplectic structures for $n>1$, i.e., $\#n\C P^2$ simply does not admit log-symplectic structures {\it bona fide} or not.

\section{Birth of singular loci}\label{sec:birth}

This section we show that under mild assumptions one can transform a symplectic structure into a log-symplectic structure with non-empty singular locus. As a consequence of this seemingly inoffensive fact we conclude that $\#m\C P^2 \# n\bar{\C P^2}$ has a log-symplectic structure with non-empty singular locus as long as $m >0 $ and  $n> 0$.

\begin{theorem}\label{theo:birth}
Let $(M^{2n},\omega)$ be a symplectic manifold and $k>0$ be an integer. If $M$ has a compact symplectic submanifold $F^{2n-2}\subset M$ with trivial normal bundle, then $M$ has a log-symplectic structure  for which the zero locus of  $\pi^n$ has $k$ components all diffeomorphic to $F\times S^1$.  
\end{theorem}

\begin{proof}
Due to the symplectic \nhood\ theorem, $F$ has a tubular \nhood\ diffeomorphic to $D^2 \times F$ endowed with the product symplectic structure, where $D^2$ is the 2-disc of radius $\e>0$. To prove the theorem it is enough to endow $D^2$ with a log-symplectic structure  whose singular locus has  $k$ components and which agrees with the standard symplectic structure near the boundary of the disc. Indeed, in this case we can consider $D^2\times F$ with the product of the log-symplectic structure on $D^2$ and the symplectic structure on $F$. Since this new structure agrees with the original symplectic structure on the boundary of the disc, we can extend it to $M$ using the original symplectic structure.

To produce the desired log-symplectic structure on $D^2$ we observe that in two dimensions every bivector is automatically Poisson, hence all we need to do is find a bivector in $D^2$ with the desired number of nondegenerate zeros. To achieve this we let $\pi \in \Gamma(\wedge^2TD)$ be the inverse of the standard symplectic structure on $D^2$ and consider the bivector $f(|x|) \pi(x)$ where $f$ is a smooth real function defined on the closed interval $[0,\e]$ which is locally constant and equal to 1 in a \nhood\ of \e, locally constant and nonvanishing in a \nhood\ of 0 and has precisely $k$  transverse zeros (see Figure \ref{fig:graph}). Then $f \pi$ is a log-symplectic structure of the desired type on $D^2$.

\begin{figure}[h!!]
\begin{center}
{\psfig{file=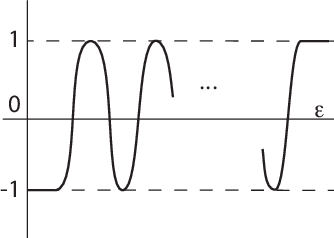, height=3.5cm,clip=}} \caption{Graph of a possible scaling function that can be used to create a singular locus with an odd number of components.} \label{fig:graph}
\end{center}
\end{figure}
\end{proof}

\begin{corollary}\label{cor:mCP2}
For any  positive integers $m,n$ the manifolds $\#m \C P^2 \#n \bar{\C P^2}$ have a log-symplectic structure whose singular locus is  diffeomorphic to $S^1 \times S^2$.
\end{corollary}
\begin{proof}
The blow-up of $\C P^2$ at a point, i.e., $\C P^2 \# \bar{\C P^2}$, has the structure of a symplectic $\C P^1$ fibration over $\C P^1$. In particular, the fibers satisfy the properties of Theorem \ref{theo:birth} and hence we can endow $\C P^2 \# \bar{\C P^2}$ with a log-symplectic structure with non-empty singular locus, say, with one component diffeomorphic to $S^1 \times S^2$. Therefore, the top power of the log-symplectic form on the symplectic locus agrees with the orientation of $\C P^2 \# \bar{\C P^2}$ at some points and disagrees in other points. By the Symplectic Blow-up Theorem \cite{MR772133}, we can blow up points in the symplectic locus and the result still has a log-symplectic structure. If we blow up points in the symplectic locus where the orientation of the log-symplectic form agrees with the orientation of $\C P^2 \# \bar{\C P^2}$, we are performing a connected sum with $\bar{\C P^2}$, while if we blow up points in the symplectic locus where the log-symplectic form gives the opposite orientation we are performing a connected sum with $\C P^2$.
\end{proof}

Notice that the manifolds obtained in Corollary \ref{cor:mCP2} have vanishing Seiberg--Witten invariants and, for $m$ and $n$ even, those manifolds do not admit almost complex structures for either choice of orientation. These are contrasts between log-symplectic and symplectic geometries since symplectic manifolds have nonzero Seiberg--Witten  invariants \cite{MR1306023} and admit almost complex structures.

As for higher dimensions, Donaldson proved that every symplectic manifold admits a Lefschetz pencil \cite{MR1802722} and hence is related to a Lefschetz fibration via the blow-up of the base locus of the pencil. Due to Theorem \ref{theo:birth} any such fibration has log-symplectic structures.

\section{Surgeries for log-symplectic manifolds}\label{sec:surgeries}

\setcounter{construction}{0}

This section we introduce surgeries which produce new log-symplectic structures out of old ones. A main feature is that in these surgeries we create new components of the singular locus hence even if the starting manifolds are symplectic the results will be only log-symplectic.

\refstepcounter{construction}\label{cons:singular fiber}
\subsection{Construction \arabic{construction}}\label{subsec:singular fiber}

This first construction produces (possibly) orientable log-symplectic manifolds out of pairs with matching data. 

\begin{buildingblock} Using the language of Theorem \ref{theo:general local form}, the local model that gives rise to the construction corresponds to the case when $Z$ has trivial normal bundle. Given a  cosymplectic manifold $(Z,\sigma,\theta)$, we let 
$$\mc{N} = (-2,2) \times Z$$
 and endow $\N$ with a log-symplectic structure for which $\{0\}\times Z$ is the singular locus, namely, we consider the 2-form 
\begin{equation}\label{eq:standard model}
\Omega = d\mathrm{log} |x| \wedge \theta + \sigma,
\end{equation}
where $|x|$ denotes the absolute value of the real number $x$.
\end{buildingblock}

\begin{ingredients} To perform this surgery we will need a (not necessarily connected) log-symplectic manifold $(M^{2n},\pi)$  together with two embeddings of a compact, connected, cosymplectic manifold $(Z,\sigma,\theta)$, $\iota_i:Z^{2n-1} \hookrightarrow M$, such that
\begin{enumerate}
\item Each $\iota_i(Z)$  lies in the symplectic locus of $\pi$ and $\iota_1(Z) \cap \iota_2(Z) =\emptyset$;
\item There is $f \in \C^{\infty}(Z)$ such that $\iota_1^*\omega = \iota_2^*\omega - df\wedge \theta = \sigma$, where $\omega$ is the symplectic form on the symplectic locus of $M$
\end{enumerate}
\end{ingredients}

\begin{thesurgery}  Since each $\iota_i(Z)$ is in the symplectic locus of $\pi$, the log-symplectic structure on $M$ gives rise to an orientation of a \nhood\ of $\iota_i(Z)$. Since $Z$ is cosymplectic, it has a natural orientation as defined by the volume form $\theta\wedge \sigma^{n-1}$.  Together the orientation on $Z$ and the (semi-local) orientation on $M$ allow us to orient the normal bundle of $Z$ and define an interior and an exterior region within the normal bundle: a vector $N \in T_{\iota_i(p)}M$ is  outward pointing if for any positive basis $\{v_1,\cdots, v_{2n-1}\}\in T_pZ$ the set $\{N,\iota_{i*}v_1,\cdots, \iota_{i*} v_{2n+1}\}$ is a positive basis for $T_pM$. Let $\hat M$ be the real oriented blow up of $M$ along $\iota_1(Z)$ and $\iota_2(Z)$, that is, $\hat M$ is diffeomorphic to the manifold obtained from $M$ by removing an open tubular \nhood\ of both copies of $Z$ and $\hat M$ has four copies of $Z$ as boundary.  At each boundary copy of $Z$, $\hat M$ lies either in the interior or the exterior side of the boundary according to the semi-local orientation. We let $\tilde M$ be the manifold obtained from $\hat M$ by identifying with each other the boundary components for which $\hat M$ lies in the interior  and similarly for the components for which $\hat M$ lies in the exterior.

\begin{figure}[h!!]
\begin{center}
{\psfig{file=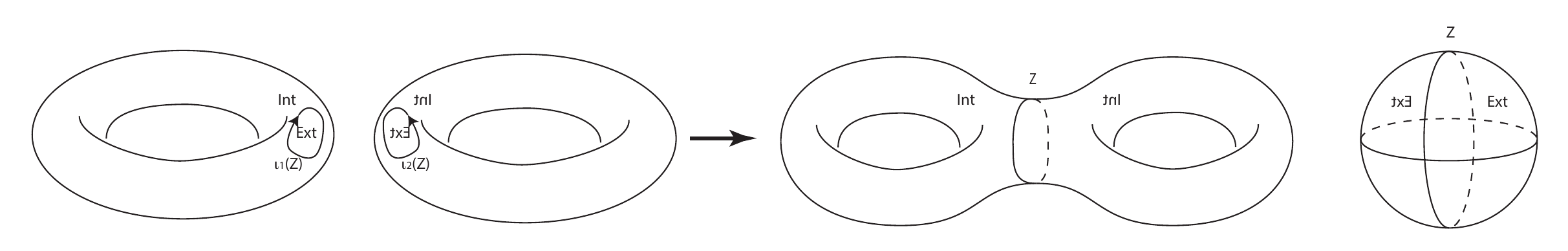, height=2.5cm,clip=}} \caption{A possible surgery on two null homologous circles lying on 2-tori. The first torus is oriented by the outward vector while the second is oriented by the inward normal vector. Interior (Int) and exterior (Ext) determined by the circles are marked in the figure with the letters inverted for different orientations. The result of the surgery is a genus two surface and a sphere.} \label{fig:graph}
\end{center}
\end{figure}

\begin{figure}[h!!]
\begin{center}
{\psfig{file=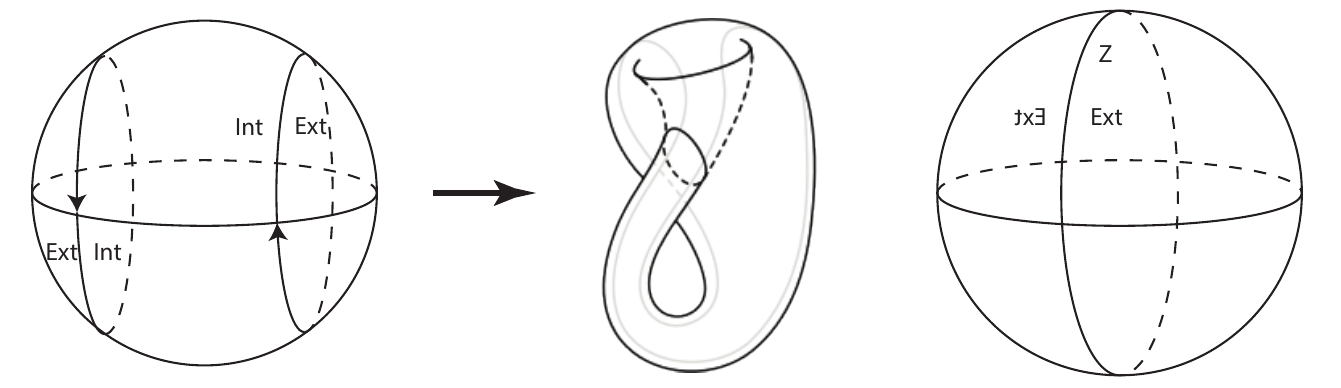, height=2.5cm,clip=}} \caption{A possible surgery on two oppositely oriented circles on a sphere yields a sphere and a Klein bottle.} \label{fig:graph}
\end{center}
\end{figure}

\begin{theorem}\label{theo:surgery1}
Let $(M,\pi)$, $(Z,\theta,\sigma)$ and $ \iota_1,\iota_2:Z\into M $ be the ingredients for the surgery and let $\hat M$ be the real oriented blow-up of $M$ along the two copies of $Z$. Then the manifold
$$\tilde M = \hat M/\sim$$
obtained by identifying the boundary components of $\hat M$ for which $\hat M$ lies in the interior (respectively exterior) of the boundary via the map $\iota_2\circ\iota_1^{-1}$ has a log-symplectic structure which agrees with the original structure on $M$ outside a \nhood\ of  two copies of  $Z = \del \hat M/\sim$ and for which $Z$ is part of the singular locus.
\end{theorem}
\begin{proof}
We have an embedding  $j_1:Z\hookrightarrow\N$,  $p\mapsto  (-1,p)$. For this embedding, $Z$ lies in the symplectic locus of the log-symplectic structure and the restriction of the symplectic form \eqref{eq:standard model} to $Z$ is just $\sigma$. Similarly, given a real function $f:Z \to \R$, for any $\e>0$  small enough we have an embedding  $j_2:Z\hookrightarrow\N$, $x \mapsto (\e e^f(p),p)$ and the restriction of the log-symplectic form to this embedding is $df \wedge \theta+\sigma$. For both embeddings, $j_1$ and $j_2$, the vector field $x \del_x$ is outward pointing with respect to the orientations induced by the symplectic and cosymplectic structures in a \nhood\ of the embeddings, that is, the cylinder
$$C = \{(x,p)\in \mc{N}: 1- \leq x \leq \e e^f(p)\}$$
contains interior points for both boundaries \wrt\ the semi-local orientations.

Hence, by Weinstein's coisotropic \nhood\ theorem, a \nhood\ of $j_1(Z)$ is symplectomorphic to a \nhood\ of $\iota_1(Z)$ and a \nhood\ of $j_2(Z)$ is symplectomorphic to a \nhood of $\iota_2(Z)$. Using these symplectomorphisms, we can glue the exterior regions of $\iota_i(Z)$ in $\hat M$ to $C$ along the boundaries and the resulting manifold has a log-symplectic structure.

We can repeat the same argumet to glue the interior regions, but now using the log-symplectic structure $-d \log |x| \wedge \theta + \sigma$ on $\mc{N}$, therefore obtaining a log-symplectic structure on $\tilde M$.
\end{proof}

\begin{remarks}
\item Even if $M$ is a symplectic manifold, and hence has a preferred orientation,  the diffeomorphism used to glue the two boundaries together does not respect these orientations hence $\tilde{M}$ does not have a preferred orientation.
\item  A common use of the theorem is when $M$ has two connected components, the maps $\iota_i$ map $Z$ to different components and there their images are separating submanifolds. In this case, $\tilde{M}$ also has two components and we will often focus our attention in one of the two, say, the one obtained by gluing the exterior regions.
\item Additive properties of the Euler characteristic imply that
$$\chi_M = \chi_{\tilde{M}}.$$
\end{remarks}

As far as examples go, Theorem \ref{theo:surgery1} leaves us with the question of how to find suitable submanifolds $Z$ to which it can be applied. Next we identify two situations in which manifolds with the desired structure appear naturally. We start with the simplest setting:

\begin{corollary}\label{theo:fiber sum}
Let $(M^{2n},\omega)$ be a log-symplectic manifold and $\iota_i:(F^{2n-2},\sigma)\into (M,\omega)$, $i=1,2$  be embeddings of a compact symplectic manifold $F$ in the symplectic locus of $M$ for which the images have trivial normal bundle. Let
$$ \tilde{M} = M\backslash(\N_1\cup \N_2)/\sim,$$
where $\sim$ indicates  the natural identification of the boundaries $\del\N_1 \cong \del \N_2$. Then $\tilde{M}$ has a log-symplectic structure which agrees with the original structure outside a tubular \nhood\ of $\del\N_i \subset \tilde{M}$. The Euler characteristic of $\tilde{M}$ is
$$\chi_{\tilde{M}} = \chi_{M}- 2 \chi_{F}.$$
\end{corollary}
\begin{proof}
Weinstein's symplectic \nhood\ theorem implies that $\iota_i(F)$ have \nhoods\ symplectomorphic to $D^2 \times F$ with the product symplectic structure. In particular the boundary of such \nhoods\ are diffeomorphic to $Z = S^1 \times F$ with the cosymplectic structure given by  the volume from of $S^1$,$\theta$, and $\sigma = \omega|_F$ and the restriction of the symplectic form to $Z$ is simply the symplectic form of $F$ pulled back to the product. Hence we  can use the theorem to conclude that $\tilde{M}$, obtained by gluing the exterior regions of $M$ \wrt\ the embeddings two embeddings of $Z$, has a log-symplectic structure.

The last claim follows from the additive properties of the Euler characteristic or by observing that the manifold obtained by gluing the interior regions (the other component of the surgery) is just $S^2 \times F$ which has Euler characteristic $2 \chi_F$. 
\end{proof}

Next we present a setting which is a little more elaborate:

\begin{corollary}\label{cor:lefschetz fibrations}
Let $p_i:M_i\into \Sigma_i$, $i = 1,2$, be symplectic Lefschetz  fibrations with the same generic fiber $(F,\sigma)$. Let $\gamma_i:\R/\lambda\Z \into \Sigma_i$  be separating loops which avoid the critical values of $p_i$, let $M_i^+$ be the exterior region of $p_i^{-1}(\gamma_i)$ and let $\gf_i:F\into F$ be the symplectomorphism of $F= p_i^{-1}(\gamma_i(0))$ obtained from symplectic parallel transport along $\gamma_i$. If there is a symplectomorphism $\vartheta:F \into F$ such that $\gf_2=\vartheta \circ \gf_1 \circ \vartheta^{-1}$, i.e., these loops have the same monodromy, then
$$M_1^+ \cup M_2^+/\sim$$
has a log-symplectic structure with singular locus $\del M_1^+ \cong \del M_2^+$.
\end{corollary}
\begin{proof}
Under the hypothesis, $Z_1 = p_1^{-1}(\gamma_1)$ is given by the quotient of $\R \times F$ by the $\Z$-action generated by
$$(x,z) \cong (x+ \lambda, \gf_1(z)).$$
The form $\theta = dx$ together with the symplectic form $\sigma \in \Omega^2(F)$ make $Z_1$  into a cosymplectic manifold.

Similarly, $Z_2 = p_2^{-1}(\gamma_2)$ is a cosymplectic manifold and the map
$$Z_1\into Z_2,\qquad (x,z) = (z, \vartheta(z))$$
is an isomorphism of cosymplectic structures as long as $\gf_2=\vartheta \circ \gf_1 \circ \vartheta^{-1}$. The result now follows from Theorem \ref{theo:surgery1}.
\end{proof}

\begin{example}[Log-symplectic structures on $\#n S^2 \times S^2$]\label{theo:nS^2}
Next we provide an explicit log-symplectic structure on $\#n S^2 \times S^2$. Our starting point is the elliptic surface $E(2k)$: the fiber sum of $2k$ copies of $\C P^2 \# 9 \bar{\C P^2}$. This is a Lefschetz  fibration $p:E(2k) \into \C P^1$ which, after appropriate identifications, has $24k$ singular fibers for which the vanishing cycles correspond to two basis  elements $\{a,b\} \in H_1(F)$ appearing in an alternating fashion, as depicted, for $E(2)$, in Figure \ref{fig:E(2)} (c.f. \cite{MR1707327}, Example 8.2.11). 

\begin{figure}[h!!]
\begin{center}
{\psfig{file=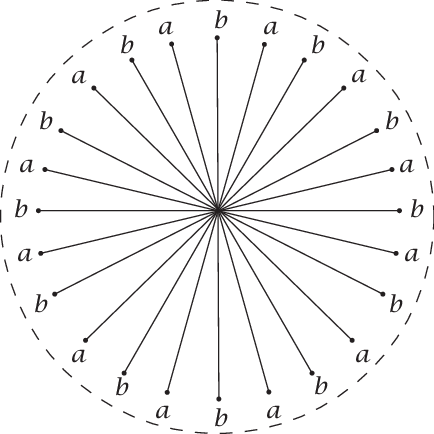, height=4.25cm,clip=}} \caption{Graphic representation of the base of $E(2)$, showing the singular values  of the projection map,  a regular value of the map for reference and paths connecting the regular value to the singular values to determine the homology class of the vanishing cycles. In this case, $a$ and $b$, the vanishing cycles which form an integral basis for $H_1(F;\Z)$.} \label{fig:E(2)}
\end{center}
\end{figure}

In order to use Construction \ref{cons:singular fiber} we consider two copies of $E(2k)$ and in both of them consider the same path, namely one whose exterior contains $n+1$ consecutive singular fibers (see Figure \ref{fig:E(n)}).

\begin{figure}[h!!]
\begin{center}
{\psfig{file=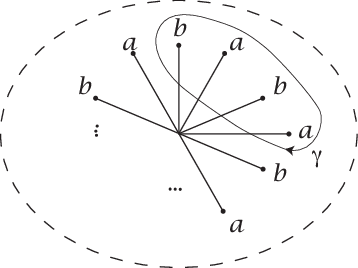, height=3.5cm,clip=}} \caption{Graphic representation of the base of $E(2k)$ and a path whose exterior contains four consecutive singular values of the projection map.} \label{fig:E(n)}
\end{center}
\end{figure}
Since we are starting with two copies of the same data we can use Corollary \ref{cor:lefschetz fibrations} to introduce a log-symplectic structure 
$$\tilde{M} =  M^+ \cup_\del M^+,$$
where  $M^+ =  p^{-1}(\overline{\Sigma^+})$ is the inverse image of the closure of the exterior points of $\gamma$. Hence  we  see that Construction \ref{cons:singular fiber} consists of taking two identical copies of $M^+\subset E(2k)$ and then glue them along the boundaries using the identity map. The resulting manifold is also known as the {\it double} of $M^+$. To precisely determine $M^+$ we observe two simple facts:
\begin{itemize}
\item  Firstly, $M^+$ admits a handlebody decomposition which contains a zero handle, $n$ 2-handles and no handles of other indices (that is $M$ is a 2-handlebody). Indeed, according to \cite{MR1707327}, Example 8.2.8, a Kirby diagram for $M$ has $n+2$ 2-handles and two 1-handles and one zero handle, but we can cancel two of 2-handles against the two 1-handles to obtain the desired handlebody decomposition);
\item Secondly, since the interior of $M^+$ is an open subset of $E(2k)$ and $E(2k)$ is spin, the intersection form on $M^+$ is even.
\end{itemize}

These facts, together with the following Proposition determine the result of the surgery:
\begin{proposition}[\cite{MR1707327}, Corollary 5.1.6]\label{prop:the double}
Let $M^4$ be a 2-handlebody with n 2-handles. Then the double of $M$ is diffeomorphic to $\#nS^2 \times S^2$ if the intersection form of $M$ is even and to $\# n \C P^2 \# n\bar{\C P^2}$ if the intersection form of $M$ is odd.
\end{proposition}

In our case, we conclude that $\tilde{M}$ is diffeomorphic to $\#n S^2 \times S^2$.
\end{example}

\end{thesurgery}

\subsection{Achiral Lefschetz fibrations}

Related to the construction of Corollary \ref{cor:lefschetz fibrations} is the notion of an achiral Lefschetz fibration:

\begin{definition}
Let $M^{2n}$ be a manifold. An {\it achiral Lefschetz fibration} on $M$ is a proper, smooth  map $p:M \into \Sigma^2$ such that the pre-image of any critical value has only one critical point and for any such pair of critical value, $y$, and critical point, $x$, there are complex coordinate systems centered at $x$ and $y$ for which $p$ takes the following form
$$p(z_1,\cdots, z_n) = z_1^2 + \cdots +z_n^2.$$
\end{definition}

Notice that in this definition we do not require $M$ or $\Sigma$ to be orientable. If they are, one can assign a sign to each critical point $x$: we demand that the complex structure on $\Sigma$ is compatible with the orientation and then we say that $x$ is {\it positive} if the complex structure on $M$ used in the definition is compatible with the orientation of $M$ and {\it negative} otherwise.

Given the construction of Corollary \ref{cor:lefschetz fibrations} and the ensuing example, one might expect that achiral Lefschetz fibrations  are related to log-symplectic structures in the same way that Lefschetz fibrations are related to symplectic structures. This is indeed the case, as we show next:

\begin{theorem}\label{theo:achiral}
Let $M^4$ and $\Sigma^2$ be compact connected manifolds and  $p:M \into \Sigma$ be an achiral Lefschetz fibration with generic fiber $F$. If $F$ is orientable and $[F]\neq 0 \in H_2(M;\R)$, then $M$ has a log-symplectic structure whose singular locus has one component and for which the fibers are symplectic submanifolds of the symplectic leaves of the Poisson structure.
\end{theorem}
\begin{proof} The proof follows closely that  of Gompf's Theorem relating Lefschetz fibrations and symplectic structures (\cite{MR1707327}, Theorem 10.2.18). Before we delve into the proof we will fix some notation. We let
\begin{itemize}
\item  $F_y= p^{-1}(y)$ for $y \in \Sigma$;
\item $\Delta$  be the set  of singular points of $p$;
\item $\Delta'$ be the set of  the singular values of $p$;
\item $\Sigma_0 = \Sigma \backslash \Delta'$ and $M_0 = p^{-1}(\Sigma\backslash \Delta')$, so that $p:M_0 \into \Sigma_0$ is a proper  fibration.
\end{itemize}

Firstly, we observe that we can assume that $p$ has connected fibers. Indeed, for achiral Lefschetz fibrations we have a short exact sequence of homotopy groups:
\begin{equation}\label{eq:sequence}
\pi_1(F)\into \pi_1(M)\stackrel{p_*}{\into} \pi_1(\Sigma)\into \pi_0(F)\into \{0\}.
\end{equation}
Since $M$ is compact, $\pi_0(F)$ is finite  and hence \eqref{eq:sequence} implies that $p_*(\pi_1(M))$has finite index in $\pi_1(\Sigma)$. Let $\tilde\Sigma$ be the cover of $\Sigma$ corresponding to the subgroup $p_*(\pi_1(M))\subset\pi_1(\Sigma)$. Then $\tilde{\Sigma}$ is compact and the map $p:M\into \Sigma$ lifts to a map $\tilde{p}:M\into \tilde{\Sigma}$. The projection $\tilde{p}$ is still an achiral Lefschetz fibration and, by construction, $\tilde{p}_*(\pi_1(M)) = \pi_1(\tilde{\Sigma})$, hence \eqref{eq:sequence} implies that the fibers of $\tilde{p}$ are connected. From now on we assume that the fibers of $p:M\into \Sigma$ are connected.

Next we deal with the general lack of orientation of the manifolds involved. Firstly,  since $F$ is orientable, we fix an orientation for $F$ for the remainder of this proof. If $p:M_0\into \Sigma_0$ was a nonorientable fibration, there would be a loop $\alpha:(I,\del I)\into M_0$ based at some point $y \in \Sigma$, where $I$ is the unit interval,  for which parallel transport (after a choice of connection) provided an orientation reversing diffeomorphism of $F_y$. In this case $p^{-1}(\alpha(I))$ would provide a chain whose boundary is $2 [F_y]$, contradicting the condition $[F]\neq 0 \in H_2(M;\R)$. Therefore $p:M_0\into \Sigma_0$ is an orientable fibration. Further this orientation induces orientations on the singular fibers of $p$ and hence we can also integrate forms over the (components of the singular) fibers. It also follows that $M$ is orientable \iff\ $\Sigma$ is.

If $\Sigma$ and $M$ are orientable, after choosing orientations, we can split the critical points of $p$ into positive and negative ones. If there are no positive or no negative points,  Gompf's Theorem (\cite{MR1707327}, Theorem 10.2.18)  implies that $M$ admits a symplectic structure and due to  Theorem \ref{theo:birth} it admits a log-symplectic structure with the desired properties. If there are positive and negative critical points, we can choose a separating loop  $\Gamma\subset  \Sigma_0$ whose interior locus contains all the negative points and whose exterior locus contains all the positive points.

If $\Sigma$ is nonorientable, we can choose a loop $\Gamma\subset \Sigma_0$ such that $\Sigma \backslash \Gamma$ is orientable and hence so is $M\backslash p^{-1}(\Gamma)$. After choosing orientations on both $\Sigma \backslash \Gamma$ and $M\backslash p^{-1}(\Gamma)$ we may homotope the loop $\Gamma$ through the negative critical  values of $p$ so that all singular points of  $p:M\backslash p^{-1}(\Gamma)\into \Sigma \backslash \Gamma$  are positive. In either case, $M\backslash p^{-1}(\Gamma)$ is an oriented manifold and $p:M\backslash p^{-1}(\Gamma) \into \Sigma \backslash \Gamma$ is a proper Lefschetz fibration possibly after changing the orientation of one of the components of  $M \backslash \Gamma$  and $\Sigma \backslash \Gamma$. From now on  we orient both $M \backslash p^{-1}(\Gamma)$  and $\Sigma \backslash \Gamma$ so that $p$ is a Lefschetz fibration there.

The next steps aim to construct a closed $2$-form $\sigma$ on $M$ which restricts to a symplectic form in every fiber.
\begin{lemma}Under the hypothesis of Theorem \ref{theo:achiral}, there is a closed 2-form $\zeta \in \Omega^2(M)$ such that
\begin{enumerate}
\item $\int_S \zeta =1 $ over any fiber and
\item if a singular fiber $F_y$ is a plumbing of two surfaces $S_1$ and $S_2$,  and we are given $s_y \in (0,1)$ then $\zeta$ can be chosen so that $\int_{S_1}\zeta = s_y$. 
\end{enumerate}
\end{lemma}
\begin{proof}
Since $[F]\neq 0 \in H_2(M;\R)$, there is a closed form  $\xi\in \Omega^2(M)$ which integrates to $1$ over the generic fibers and hence over all fibers. Next we need to argue that one can change $\xi$ so that the property 2. holds. In this case, with the orientations chosen before on $M \backslash p^{-1}(\Gamma)$, the intersection number of $S_1$ and $S_2$ is 1. Since $S_1 \cup S_2$ is homologous to a regular fiber, $\int_{S_1} \xi+ \int_{S_2} \xi =1$.   If $\int_{S_1}\xi = r$, let $\psi$ be a form with support in a \nhood\ of $S_2\subset M\backslash p^{-1}(\Gamma)$ which represents the Poincar\'e dual of $S_2$ and consider $\xi' = \xi  + (-r +s_y)\psi$. Then for any closed surface $S'$ inside another fiber $F_{y'}$,  $\int _{F_{y'}}\psi=0 $ as $F_{y'}$ does not intersect $S_2$ hence $\int_{S'}\xi' = \int_{S'}\xi$.  On the other hand
$$\int_{S_1} \xi' = \int_{S_1} \xi \pm (-r+s_y)\int_{S_1} \psi = s_y.$$
That is, after a change in $\xi$, we  found another closed form for which the claim holds at a specific singular fiber. Since there are only finitely many such fibers, we can repeat the process for each of them to obtain the desired $\zeta$. 
\end{proof}

\begin{lemma}
Under the hypothesis of Theorem \ref{theo:achiral}, there is a finite good open  cover $\U$ of $\Sigma$ such that for each $U_\alpha \in \U$, there is a closed form $\eta_\alpha \in \Omega^2(p^{-1}(U_\alpha))$ which is symplectic on the fibers of $p$.   
\end{lemma}
\begin{proof}
Since  $p:M\backslash p^{-1}(\Gamma) \into \Sigma \backslash \Gamma$ is a proper Lefschetz fibration  for which $[F]\neq 0$,  it follows from Gompf's Theorem (\cite{MR1707327}, Theorem 10.2.18) that each component of $M\backslash p^{-1}(\Gamma)$ has a symplectic form for which the fibers are symplectic of area 1. Letting $\N$ be a small tubular \nhood\ of $\Gamma$ without singular values of $p$, it follows that, $p:p^{-1}(\N) \into \N$ is a proper Lefschetz fibration for which $[F]\neq 0$ hence we can also apply Gompf's result here to conclude that there is a symplectic form $\omega_0$ on $p^{-1}(\N)$ for which the fibers are symplectic of area 1. Finally, let $\U$ be a finite good refinement of the cover $\{\Sigma\backslash \Gamma,\N\}$ and for each $U_\alpha \in \U$  let $\eta_\alpha$ be the restriction of one of the symplectic forms above to $p^{-1}(U_\alpha)$.  
\end{proof}

\begin{lemma}
Under the hypothesis of Theorem \ref{theo:achiral}, there is a closed 2-form $\sigma \in \Omega^2(M)$ such that $\sigma|_{F_y}$ is a symplectic form for every fiber $F_y$.  
\end{lemma}
\begin{proof}
The forms  $\eta_\alpha$ and $\zeta$ are cohomologous on  $p^{-1}(U_\alpha)$  (in the case of singular fibers, one must choose the value $s_y$ so that the integrals of these forms over each cycle agree). Therefore  there are $\theta_i \in \Omega^1(p^{-1}(U_i))$ such that $\eta_i = \zeta + d\theta_i$. Let $\{\kappa_\alpha\}$ be a partitition of unity subbordinate to the cover $\U$ and consider the form
$$\sigma = \zeta + d\sum_i p^*(\kappa_i) \theta_i.$$
Since $p^*\kappa|_{F_y}$ is a constant, we have that 
$$\sigma|_{F_y} = \zeta|_{F_y} + \sum_i p^*(\kappa_i) d\theta_i|_{F_y} =  \sum_i p^*(\kappa_i)(\zeta|_{F_y}+ d\theta_i|_{F_y}) = \sum_i p^*(\kappa_i) \eta_i|_{F_y}.$$
Since each $\eta_i|_{F_y}$ is a symplectic form and all of them determine the same orientation $\sum_i p^*(\kappa_i) \eta_i|_{F_y}$ is a symplectic form on $F_y$.
\end{proof}

\noindent
{\it End of proof of Theorem \ref{theo:achiral}}. 
Finally, to obtain the log-symplectic structure on $M$, observe that since $\Sigma\backslash \Gamma$ is oriented, $\Gamma$ is a real divisor on $\Sigma$ representing $\wedge^2T\Sigma$, i.e., there is a section $\pi \in \Gamma(\wedge^2(T\Sigma))$ which has $\Gamma$ as its (transverse) zero locus. Since $\Sigma$ is two-dimensional, $\pi$ is  a log-symplectic structure whose singular locus is $\Gamma$. We further choose $\pi$ so that it agrees with the orientation of $\Sigma\backslash \Gamma$. Inverting $\pi$ we obtain  $\omega_\Sigma\in \Omega^2(\Sigma)$ with a log-singularity at $\Gamma$. Then the standard argument shows that $\omega = p^*\omega_\Sigma +\e \sigma$ is a log-symplectic structure for $\e$ small enough. Indeed, away from critical points of $p$, $\omega_\Sigma$ dominates $\sigma$ and hence determines a log-symplectic structure on the complement of a small \nhood\ of $\Delta$. In particular, $\omega$ determines an orientation on its symplectic locus. With respect to this orientation on $M\backslash p^{-1}(\Gamma)$ and the orientation determined by $\pi$ on $\Sigma\backslash \Gamma$, all singular points are positive and the argument from \cite{MR1707327}, Exercise 10.2.21, shows that $\omega$ is symplectic on $M\backslash p^{-1}(\Gamma)$.
\end{proof}

\noindent
{\it Remarks} ({\it The condition  $[F]\neq 0$}). 
\begin{itemize}
\item If  the genus of the fiber is different from 1, then $\mathrm{ker}(p)$ defines a line bundle over $M \backslash \Delta$ and this line bundle extends to the singular locus. Letting $c_1$ be the first Chern class of this bundle, naturality implies that $c_1|_F$ is just the Euler class of the fiber and hence, if the genus of the fiber is not 1, $c_1$ evaluates nonzero on $[F]$, showing that $[F] \neq 0$;
\item If an achiral Lefschetz fibration over an oriented surface has a section, then any fiber represents a nontrivial class since it has nontrivial intersection with the section.
\item $S^4$ admits an achiral Lefschetz fibration. Since $H^2(S^4) = \{0\}$,  the fibers are homologically trivial and hence are tori and there is no section. Further, $S^4$ also does not admit log-symplectic structures due to Theorem \ref{marcut-Osorno-Torres}. Therefore the condition $[F]\neq 0$ can not be removed from Theorem \ref{theo:achiral}.
\end{itemize}
 
 \begin{remark}
The proof above is also very similar to the one given by Baykur \cite{MR2253445} relating Achiral Lefschetz fibrations to folded symplectic structures, as both ours and Baykur's proof follow Gompf's original proof closely (\cite{MR1707327}, Theorem 10.2.18). The main differences between the proofs regard the treatment of the singular locus, as folded and log-symplectic structures have different types of singular behaviour and Gompf did not have to deal with either of them. Further, since log-symplectic structures can always be deformed into folded ones, our proof is a little more general Baykur's.
\end{remark} 
 
 Achiral Lefschetz fibrations have been studied by Etnyre and Fuller \cite{MR2219214} and are present in several four-manifolds:
 
\begin{theorem}[Etnyre--Fuller \cite{MR2219214}]
Let $X$ be a smooth, closed, oriented 4-manifold. Then there exists a framed circle in $X$ such that the manifold obtained by surgery along that circle admits an achiral Lefschetz fibration with section and whose base is  $S^2$. Further, if $M$ is simply-connected then we can arrange so that both $M \# (S^2 \times S^2)$ and $M\#\C P^2 \# \bar{\C P^2}$ arise as such surgery and hence both $M \# (S^2 \times S^2)$ and $M\#\C P^2 \# \bar{\C P^2}$  admit achiral Lefschetz fibrations   with a section over $S^2$.  
\end{theorem}

Combining Theorem \ref{theo:achiral} with Etnyre--Fuller's Theorem we get:

\begin{theorem}\label{theo:corollary from Etnyre--Fuller}
Let $M$ be a simply connected compact four-manifold. Then both $M \# (S^2 \times S^2)$ and $M\#\C P^2 \# \bar{\C P^2}$ admit bona fide log-symplectic structures.  
\end{theorem}

\refstepcounter{construction}\label{cons:twisted blow up}
\subsection{Construction \arabic{construction}}

The second construction is a nonorientable version of the first which produces proper log-symplectic manifolds.

\begin{buildingblock}
Given a symplectic manifold $(F,\sigma)$ a symplectomorphism $\gf:F \into F$ and $\lambda> 0$ we form the quotient of the log-symplectic manifold
$$\N = (-2,2)\times \R \times F;\qquad \Omega|_{(x,y,p)} = d\log|x| \wedge dy + \sigma$$
by the $\Z$-action generated by $(x,y,p) \sim (-x,y +\lambda,\gf(p))$:
$$\N_\gf = \N/\Z.$$
Then $\N_\gf$ is a log-symplectic manifold with singular locus $Z = \{0\}\times \R \times F/\Z$. Notice that  $\mc{N}_\gf\backslash Z$ is a fiber bundle over $\R_+$, with projection map induced  by the invariant map $\pi_1(x,y,p) = |x|$ defined on $(-2,2)\times \R \times F$. Then $Z' =\pi_1^{-1}(1)$ is a coisotropic submanifold of the symplectic locus given by
\begin{equation}\label{eq:Z'}
Z' = \R\times F/\Z\qquad (y,p) \sim (y + 2\lambda, \gf^2(p)).
\end{equation}
\end{buildingblock}

\begin{ingredients} We will need a proper cosymplectic manifold $(Z',\theta,\sigma)$  with symplectic fiber $F$ for which the monodromy map is the square of a symplectomorphism $\gf:F \into F$, i.e., $Z'$ is given by \eqref{eq:Z'}. We will need further a log-symplectic manifold $(M,\pi)$ and  a separating embedding $\iota:Z' \hookrightarrow M$ in the symplectic locus of $M$ such that $\iota^*\omega = \sigma$, where $\omega$ is the induced symplectic structure on $M$.
\end{ingredients}

\begin{thesurgery}
The surgery follows the same lines of Construction \ref{cons:singular fiber}: Since $\iota(Z)$ is separating, it defines an exterior and an interior region of $M$. Let $M^+$ be the closure of the exterior. Then there  is a $\Z_2$-action on  $Z' = \del M^+$, namely, in terms of \eqref{eq:Z'}, the action is generated by the map
$$(y,p) \mapsto (y+\lambda,\gf(p))$$
and the orbits of this action form an equivalence relation on $\del M^+$ which allows us  to form the space
$$\tilde{M} = M^+/\sim$$
obtained by taking the quotient of $\del M^+$ by this equivalence relation.

\begin{theorem}\label{theo:twisted blow up}
Let $(M,\pi)$, $(Z',\theta,\sigma)$ and $ \iota:Z'\into M $ be the ingredients for the surgery and let $M^+$ be the closure of the exterior region defined by $\iota(Z')$. Then the manifold
$$\tilde{M} = M^+/\sim$$
obtained by taking the quotient of $\del M^+$ by the $\Z_2$-action has a log-symplectic structure which agrees with the original structure on $M^+$ outside a \nhood\ of  $Z = \del M^+/\sim$ and for which $Z$ is part of the singular locus.\end{theorem}

The proof of this theorem is completely analogous to that of Theorem \ref{theo:surgery1}. A particular case of this surgery has a geometric interpretation.

\begin{corollary}[Real blow-up]\label{cor:real blow up}
Let $(M^{2n},\omega)$ be a log-symplectic manifold and let $F^{2n-2}\subset M$ be a symplectic submanifold which does not intersect the singular locus and has trivial normal bundle. Then the real blow-up of $M$ along $F$ has a log-symplectic structure for which the exceptional divisor is a component of the singular locus.
\end{corollary}
\begin{proof}
Just as in Corollary \ref{theo:fiber sum}, the requirement that $F$ has trivial normal bundle implies that a \nhood\ of $F$ is symplectomorphic to $D^2 \times F$ and hence we obtain an embedding of the proper cosymplectic manifold $S^1 \times F$ into $M$. The monodromy of this cosymplectic manifold is the identity map which is obviously the square of a symplectomorphism. Now, the local model is based on using $\gf = \Id$, that is $\N_F = \mathbb{M}\times F$, where $\mathbb{M}$ is the M\"obius band and the effect of the surgery is that we remove a \nhood\ of $F$ (which is diffeomorphic to $D^2 \times F$) and glue back $\mathbb{M}\times F$. This is precisely the underlying surgery of the real blow-up of $F$.
\end{proof}

\end{thesurgery}

\section{Reversing the surgeries}\label{sec:reversing}

Last section we managed to produce several examples of log-symplectic manifolds out of symplectic manifolds. One might rightfully expect that there are more examples of such structures: for one thing the Stiefel--Whitney class either vanished (first construction) or corresponded to the generator of  $H^1(S^1;\Z_2)$ (second construction),  therefore leaving out a number of possibilities. On the other hand, if we assume that $M$ is orientable or, in the nonorientable case, take the orientable double cover, then any singular locus automatically is associated to the zero Stiefel-Whitney class and hence it has neighboorhood  diffeomorphic, as a Poisson, manifold to the building blocks used in Construction \ref{cons:singular fiber}. Next we show that in four dimensions any log-symplectic structure is created out of our surgeries and hence can be cut up and filled into a collection of compact symplectic manifolds.

\begin{theorem}\label{theo:4d}
Let $(M^4,\pi)$ be a compact, orientable, log-symplectic manifold with singular locus $Z$. Then each unoriented component of $M\backslash Z$ can be compactified as a symplectic manifold.
\end{theorem}
\begin{proof}
According to Theorem \ref{theo:general local form}, a \nhood\ of each connected component of $Z$ is equivalent to the building block of  Construction \ref{cons:singular fiber} and hence we have two copies of $Z$ in such \nhood\ (one on either side of the singular locus) as a coisotropic submanifold. To reverse the surgery, one would need to prove that such coisotropic submanifold appears as the boundary of the (interior of) a symplectic manifold. But in four dimensions any cosymplectic manifold is automatically a taut foliation and hence the conclusion follows from the following theorem:
\begin{theorem}[\cite{MR2388043}, Theorem 41.3.1]
Let $Z$ be a closed 3-manifold and $\mc{F} \subset TZ$ be a smooth taut foliation. Let $\sigma \in \Omega^2(Z)$ be the closed form which is positive on the leaves of $\mc{F}$. Then there is a closed symplectic manifold $(X,\omega)$ containing $Z$ as a separating submanifold such that  $\omega|Z = \sigma$.
\end{theorem}
\end{proof}

\bibliographystyle{hyperamsplain}
\bibliography{references}

\providecommand{\bysame}{\leavevmode\hbox to3em{\hrulefill}\thinspace}
\providecommand{\MR}{\relax\ifhmode\unskip\space\fi MR }
\providecommand{\MRhref}[2]{%
  \href{http://www.ams.org/mathscinet-getitem?mr=#1}{#2}
}
\providecommand{\href}[2]{#2}
\begin{thebibliography}{10}

\bibitem{MR2253445}
R.~{\.I}. Baykur, \emph{K\"ahler decomposition of 4-manifolds},
  \href{http://dx.doi.org/10.2140/agt.2006.6.1239}{Algebr. Geom. Topol.
  \textbf{6} (2006)}, 1239--1265.

\bibitem{MR2670164}
A.~Cannas~da Silva, \emph{Fold-forms for four-folds}, J. Symplectic Geom.
  \textbf{8} (2010), no.~2, 189--203.

\bibitem{MR1748286}
A.~Cannas~da Silva, V.~Guillemin, and C.~Woodward, \emph{On the unfolding of
  folded symplectic structures},
  \href{http://dx.doi.org/10.4310/MRL.2000.v7.n1.a4}{Math. Res. Lett.
  \textbf{7} (2000)}, no.~1, 35--53.

\bibitem{MR1802722}
S.~K. Donaldson, \emph{Lefschetz pencils on symplectic manifolds}, J.
  Differential Geom. \textbf{53} (1999), no.~2, 205--236.

\bibitem{MR2219214}
J.~B. Etnyre and T.~Fuller, \emph{Realizing 4-manifolds as achiral {L}efschetz
  fibrations}, \href{http://dx.doi.org/10.1155/IMRN/2006/70272}{Int. Math. Res.
  Not. (2006)}, Art. ID 70272, 21.

\bibitem{FMM13}
P.~Frejlich, D.~Martinez-Torres, and E.~Miranda, \emph{Constructing global
  $b$-symplectic examples via symplectic topology}, 2013. Work in progress.

\bibitem{MR1707327}
R.~E. Gompf and A.~I. Stipsicz, \emph{{$4$}-manifolds and {K}irby calculus},
  Graduate Studies in Mathematics, vol.~20, American Mathematical Society,
  Providence, RI, 1999.

\bibitem{Gualtieri-Li-2012}
M.~Gualtieri and S.~Li, \emph{Symplectic groupoids of log symplectic
  manifolds}, 2012. \href{http://arxiv.org/abs/arXiv:1206.3674}{{\tt
  arXiv:1206.3674}}.

\bibitem{MR2861781}
V.~Guillemin, E.~Miranda, and A.~R. Pires, \emph{Codimension one symplectic
  foliations and regular {P}oisson structures},
  \href{http://dx.doi.org/10.1007/s00574-011-0031-6}{Bull. Braz. Math. Soc.
  (N.S.) \textbf{42} (2011)}, no.~4, 607--623.

\bibitem{Guilleminetal2012}
\bysame, \emph{Symplectic and Poisson geometry on $b$-manifolds}, 2012.
  \href{http://arxiv.org/abs/arXiv:1206.2020}{{\tt arXiv:1206.2020}}.

\bibitem{Guilleminetal2015}
V.~Guillemin, E.~Miranda, and J.~Weitsman, \emph{Desingularizing
  $b^m$-symplectic structures}, 2015.
  \href{http://arxiv.org/abs/arXiv:1512.05303}{{\tt arXiv:1512.05303}}.

\bibitem{MR2388043}
P.~Kronheimer and T.~Mrowka,
  \href{http://dx.doi.org/10.1017/CBO9780511543111}{\emph{Monopoles and
  three-manifolds}}, New Mathematical Monographs, vol.~10, Cambridge University
  Press, Cambridge, 2007.

\bibitem{MR2481690}
H.~Li, \emph{Topology of co-symplectic/co-{K}\"ahler manifolds},
  \href{http://dx.doi.org/10.4310/AJM.2008.v12.n4.a7}{Asian J. Math.
  \textbf{12} (2008)}, no.~4, 527--543.

\bibitem{marcut-osorno}
I.~Marcut and B.~Osorno-Torres, \emph{On cohomological obstructions for the
  existence of log-symplectic structures}, 2013.
  \href{http://arxiv.org/abs/arXiv:1303.6246}{{\tt arXiv:1303.6246}}.

\bibitem{MR772133}
D.~McDuff, \emph{Examples of simply-connected symplectic non-{K}\"ahlerian
  manifolds}, J. Differential Geom. \textbf{20} (1984), no.~1, 267--277.

\bibitem{osorno-thesis}
B.~Osorno-Torres, \emph{Ph.D. thesis. Utrecht University}, 2013. Work in
  progress.

\bibitem{MR1959058}
O.~Radko, \emph{A classification of topologically stable {P}oisson structures
  on a compact oriented surface}, J. Symplectic Geom. \textbf{1} (2002), no.~3,
  523--542.

\bibitem{MR1306023}
C.~H. Taubes, \emph{The {S}eiberg-{W}itten invariants and symplectic forms},
  Math. Res. Lett. \textbf{1} (1994), no.~6, 809--822.

\end{thebibliography}

\end{document}